\documentclass[12pt,reqno]{amsart}
\usepackage{amssymb}
\usepackage{amsmath}
\usepackage{amsthm}
\usepackage{mathptmx}
\usepackage{color}
\usepackage[pdftex]{graphicx}
\newcommand{\R}{\mathbb R}

\newcommand{\Z}{\mathbb Z}

\newcommand{\car}{\chi_{_{0,\epsilon,b}}}
\newcommand{\carn}{\chi_{_{n,\epsilon,b}}}

\newtheorem{theorem}{Theorem}[section]

\newtheorem{remark}[theorem]{Remark}

\newtheorem{corollary}[theorem]{Corollary}

\newtheorem*{TA}{Theorem A}
\newtheorem*{TB}{Theorem B}
\newtheorem*{TC}{Theorem C}
\newtheorem*{TD}{Theorem D}
\numberwithin{equation}{section}
\begin{document}
\title[Regularity and decay of the k-generalized KdV equations]{On the propagation of regularity and decay of solutions to the $k$-generalized Korteweg-de Vries equation}
%%%%%%%%%%%%%%%%%%
\author{Pedro Isaza}
\address[P. Isaza]{Departamento  de Matem\'aticas\\
Universidad Nacional de Colombia\\ A. A. 3840, Medellin\\Colombia}
\email{pisaza@unal.edu.co}
%\thanks{The third author is supported  by  the NSF  DMS-1101499}
%\thanks{}
%%%%%%%%%%%%%%%%%%%%%%
\author{Felipe Linares}
\address[F. Linares]{IMPA\\
Instituto Matem\'atica Pura e Aplicada\\
Estrada Dona Castorina 110\\
22460-320, Rio de Janeiro, RJ\\Brazil}
\email{linares@impa.br}

\author{Gustavo Ponce}
\address[G. Ponce]{Department  of Mathematics\\
University of California\\
Santa Barbara, CA 93106\\
USA.}
\email{ponce@math.ucsb.edu}
%\thanks{}
%%%%%%%%%%%%%%%%%%%%%%

%\date{7/17/14}
%\dedicatory{}
%%%%%%%%%%%%%%
\begin{abstract}
We study special regularity and decay properties of solutions to the IVP associated to the $k$-generalized KdV equations.
In particular,  for datum $u_0\in H^{3/4^+}(\R)$ whose restriction belongs to $H^l((b,\infty))$ for some $l\in\Z^+$ and $b\in \R$ we prove that the restriction of the corresponding solution 
$u(\cdot,t)$ belongs to $H^l((\beta,\infty))$ for any $\beta \in \R$ and any $t\in (0,T)$. Thus, this type of regularity propagates with infinite speed to its left as time evolves.
\end{abstract}

\maketitle
\section{Introduction}

In this work we shall deduce some special properties of solutions to the $k$-generalized Korteweg-de Vries (k-gKdV) equations. These  properties are concerned with the propagation of regularity and decay of their solutions on the line. Thus, we shall consider the initial value problem (IVP) 
\begin{equation}\label{gkdv}
\begin{cases}
\partial_t u+\partial_x^3u+u^k\,\partial_x u=0, \hskip10pt x,t\in\R, \;\;k\in \Z^{+},\\
u(x,0)=u_0(x).
\end{cases}
\end{equation}

 The cases $k=1$ and $k=2$ in \eqref{gkdv} correspond to the KdV and modified KdV (mKdV) equations, respectively. Initially they arise as a model in nonlinear wave propagation in a shallow channel and later as models in several other physical phenomena  (see \cite{Miur} and references therein). Also they have been shown to be completely integrable. In particular, they possess infinitely many conservation laws. In the case of powers $k=3,4, \dots$ the solutions of the equation in \eqref{gkdv} satisfy just three conservation laws. 

  Following \cite{Kato83} it is said that the IVP \eqref{gkdv} is locally well-posed (LWP) in the function space $X$ if given any datum $u_0\in X$ there exist $T>0$ and a unique solution
$$
u\in C([-T,T] :X)\cap\dots.
$$
of the IVP \eqref{gkdv} with the map data-solution, $u_0\to u$, being continuous. In the cases where $T$ can be taken arbitrarily large, one says that the problem is globally well-posed (GWP). In both cases, the solution flow defines a dynamical system on $X$. 

For the IVP \eqref{gkdv} the natural function space is the classical Sobolev family 
\begin{equation}
\label{sob}
 H^s(\R)=(1-\partial_x^2)^{-s/2}L^2(\R),\;\;\;\;\;s\in\R.
\end{equation}
 
 The problem of finding the minimal regularity on the data $u_0$ in the scale in \eqref{sob} which guarantees that the IVP \eqref{gkdv} is well posed has been studied extensively. After the  works of \cite{SaTe},
\cite{BoSm} and \cite{Kato83} gathering the local and the global well-posedness results, in \cite{KPV93}, \cite{Bo}, \cite{KPV96}, \cite{ChCoTa}, \cite{CKSTT}, \cite{Ta}, \cite{GuPaSi},  \cite{Guo} and \cite{Ki} one has:

The IVP \eqref{gkdv} with $k=1$ is globally well-posed in $H^s(\R)$ for $s\geq -3/4$   (see \cite{Bo} (GWP for $s=0$), \cite{KPV96} (LWP for $s>-3/4$), \cite{ChCoTa} (LWP for $s\geq-3/4$), \cite{CKSTT} (GWP for $s>-3/4$) , and \cite{Guo}, \cite{Ki}  (GWP for $s\geq -3/4$)).
 
The IVP \eqref{gkdv} with $k=2$ is globally well-posed in $H^s(\R)$ for $s\geq 1/4$   (see \cite{KPV93} (LWP for $s\geq 1/4$), \cite{CKSTT} (GWP for $s>1/4$) and 
\cite{Guo}, \cite{Ki} (GWP for $s\geq 1/4$)).

The IVP \eqref{gkdv} with $k=3$ is locally well-posed in $H^s(\R)$ for $s\geq -1/6$ (see \cite{Gu} for LWP for $s>-1/6$ and \cite{Ta} for LWP for $s=-1/6$)  and GWP  in $H^s(\R)$ for $s>-1/42$ (\cite{GuPaSi}).

The IVP \eqref{gkdv} with $k\geq 4$ is locally well-posed in $H^s(\R)$ for $s\geq (k-4)/2k$ (see \cite{KPV93}).

In \cite{MaMe} it was shown that a local solution of the IVP \eqref{gkdv} with $k=4$ corresponding to smooth initial data can develop singularities in finite time.

 The well-posedness of the IVP \eqref{gkdv} has been also studied  in the weighted Sobolev spaces
\begin{equation}
\label{wei-sob}
Z_{s,r}= H^s(\R)\cap L^2(|x|^rdx),\;\;\;\;\;s,\,r\in\R.
\end{equation}

 In \cite{Kato83} it was established the LWP of the IVP \eqref{gkdv} in $Z_{s,r}$ with $r$ an even integer and $s\geq r$. In particular, this implies the well-posedness 
 of the IVP \eqref{gkdv} in the Schwartz space $\mathcal S(\R)$. The proof of this result is based in the commutative relation of the operators 
$$
\mathcal L=\partial_t+\partial_x^3\;\;\;\;\text { and }\;\;\;\Gamma =x-3t\partial_x^2.
$$ 

 Roughly, in \cite{Na} and \cite{FoLiPo} the above result was extended to the case $Z_{s,r}$ with $r>0$, $s\geq r$ and $s\geq \max\{s_k;0\}$ with $ s_1=-3/4, \,s_2=1/4$ 
 and  $s_k=(k-4)/2k \;$ for $\,k\geq 3$. 

 In \cite{ILP1} it was shown that the hypothesis $s\geq r$ is necessary. More precisely, it was established in \cite{ILP1}  that if $u\in C([-T,T]:H^{s}(\R)) \cap\dots$ is a solution of the IVP \eqref{gkdv} 
 with $s\geq \max\{s_k;0\}$ (with $s_k$ as above) and there exist two different times $t_1,\,t_2\in [-T,T]$  such that
$ \, |x|^{\alpha}u(\cdot,t_j)\in L^2(\R)$ for $j=1,2$ with $2\alpha>s$, then $u\in C([-T,T]:H^{2\alpha}(\R))$. 

Thus, if $u_0\in Z_{s,r}$ with $r>s$, then at time $\,t\neq 0$  the solution $u(\cdot,t)$ stays only in  $Z_{s,s}$, i.e. the extra decay $r-s$ is not preserved by the solution flow. As a consequence of our results obtained here we shall see that this extra-decay, no preserved by the flow, is transformed into extra regularity in a precise manner (see Corollary \ref{cor7}).

 The above results are concerned with regularity and decay properties of solutions of the IVP \eqref{gkdv} in symmetric spaces.

 For asymmetric weighted spaces, for which the outcome 
 has to be restricted to forward times $t>0$, one has the following result found in \cite{Kato83} for the KdV equation, $k=1$ in \eqref{gkdv}, in the space $L^2(e^{\beta x}dx),\,\beta>0$. In \cite{Kato83} it was shown that the persistence property holds for $L^2$-solutions in $L^2(e^{\beta x}dx),\,\beta>0$, for $t>0$. Moreover, formally in this space the operator $\partial_t+\partial_x^3$  becomes $\partial_t+(\partial_x-\beta)^3$ so the solutions of the equation exhibit a parabolic behavior. More precisely, the following result for the KdV equation was proven in \cite{Kato83} (Theorem 11.1 and Theorem 12.1).

   \begin{TA} \label{Ka}
  Let $ u\in C([0,\infty)\,:\,H^2(\R))$ be a solution of the IVP \eqref{gkdv} with $k=1$. If 
  \begin{equation}
  \label{1.3}
  u_0\in H^2(\R)\cap L^2(e^{\beta x}dx),\;\;\;\;\text{for some}\;\;\beta>0,
 \end{equation}
  then
  \begin{equation}
  \label{1.4} 
e^{2\beta x}u\in   C([0,\infty)\,:\,L^2(\R))\cap C((0,\infty)\,:\,H^{\infty}(\R)),
\end{equation}
with
\begin{equation}
\label{1.5}
\|u(t)\|_2=\|u_0\|_2,\;\;\;\;\;\;\;\|u(t)-u_0\|_{-3,2}\leq Kt,\;\;\;\;\;\;t>0,
\end{equation}
and
\begin{equation}
\label{1.6}
\|e^{\beta x}u(t)\|_2\leq e^{Kt}\,\| e^{\beta x}u_0\|_2,\;\;\;\;\;t>0,
\end{equation}
where $K=K(\beta,\|u_0\|_2)$. Moreover, the map data-solution $u_0\to u$ is continuous from $L^2(\R)\cap  L^2(e^{\beta x}dx)$ to  $C([0,T]\,:\,L^2(e^{\beta x}dx))$, for any $T>0$.
 \end{TA}

 It is easy to see that the result of Theorem A extends to solutions of the IVP \eqref{gkdv} for any $k\in\Z^+$ in their positive time interval of existence  $\,[0,T]$. 

 The problem of the uniqueness of solutions in polynomial weighted $\,L^2(\R)$-spaces was considered in \cite{KrFa} 
and \cite{GiTs}. In particular, in \cite{GiTs} uniqueness and \it a priori \rm estimates, for positive times, were deduced for solutions of the KdV equation ($k=1$ in \eqref{gkdv})  with data $u_0\in L^2((1+x_{+})^{3/4}dx)$ and for the mKdV equation ($k=2$ in \eqref{gkdv})  for data $u_0\in L^2((1+x_{+})^{1/4}dx)$.

 In view of the results in Theorem A is natural to ask what is the strongest weighted space where persistence of the solutions of \eqref{gkdv} holds. The following uniqueness result obtained in \cite{EKPV06} gives an upper bound of this weight.  

\begin{TB} \label{EKPV06}
 There exists $c_0=c_0(k)>0$ such that for any pair  
 $$
u_1,\,u_2\in C([0,1]:H^4(\R)\cap L^2(|x|^2dx))
$$
of  solutions of  \eqref{gkdv}, if 
\begin{equation}
 \label{3:2}
 u_1(\cdot,0)-u_2(\cdot,0),\,\;\, u_1(\cdot,1)-u_2(\cdot,1)\in L^2(e^{c_0x_{+}^{3/2}}dx),
\end{equation}  then $u_1\equiv u_2$.
  \end{TB}

Thus, by taking $u_2\equiv 0$ one gets  a restriction on the possible decay 
of a non-trivial solution of 
\eqref{gkdv}
at two different times. It was established in \cite{ILP1} that this result is optimal. More precisely, the following was proven in \cite{ILP1}:

\begin{TC} \label{C}
Let $a_0$ be a positive constant. For any given data
\begin{equation}
\label{0.1}
u_0\in H^{s}(\R)\cap L^2(e^{a_0x_{+}^{3/2}}dx),\;\;\;\;s>s_k \;\;\;\text{defined above}
\end{equation}
the unique solution of the IVP \eqref{gkdv} satisfies that for any $T>0$
\begin{equation}
\label{0.2}
\sup_{t\in [0,T]}\,\int_{-\infty}^{\infty} e^{a(t)x_{+}^{3/2}}|u(x,t)|^2dx  \leq c^*
\end{equation}
with $\;=c^*( a_0; \|u_0\|_2; \|e^{a_0x_{+}^{3/2}/2}u_0\|_2; T; k) $ and 
\begin{equation}
 \label{0.3}
a(t)=\frac{a_0}{(1+27 a_0^2t/4)^{1/2}}.
 \end{equation}

\end{TC}

  Before stating our results we need  to define the class of solutions to the IVP \eqref{gkdv} to which it applies. Thus, we shall rely on the following 
  well-posedness result which is a consequence of the arguments deduced in \cite{KPV93}: 

 \begin{TD} \label{TD}
If $u_0\in H^{{3/4}^{+}}(\R)$, then there exist $T=T(\|u_0\|_{{3/4}^{+}, 2}; k)>0$ and a unique solution of the IVP \eqref{gkdv}
such that
\begin{equation}\label{notes-2}
\begin{split}
{\rm (i)}\hskip25pt & u\in C([-T,T] : H^{{3/4}^{+}}(\R)),\\
{\rm(ii)}\hskip25pt & \partial_x u\in L^4([-T,T]: L^{\infty}(\R)), \hskip10pt {\text (Strichartz)},\\
{\rm(iii)}\hskip25pt &\underset{x}{\sup}\int_{-T}^{T} |J^r\partial_x u(x,t)|^2\,dt<\infty \text{\hskip10pt for \hskip10pt} r\in[0,{3/4}^{+}],\\
{\rm(iv)}\hskip25pt & \int_{-\infty}^{\infty}\; \sup_{-T\leq t\leq T}|u(x,t)|^2\,dx < \infty,
\end{split}
\end{equation}
with $\;J=(1-\partial_x^2)^{1/2}$. Moreover, the map data-solution, $\,u_0\to \,u(x,t)$ is locally continuos (smooth) 
from $\,H^{3/4+}(\R)$ into the class defined in \eqref{notes-2}.
\end{TD} 
 
 As it was already commented, in the cases $k=1,2,3$ the value of  $T$ in Theorem D can be taken arbitrarily large. Also we shall remark 
that for powers $k=2,3,..$ it suffices to assume that $u_0\in H^{3/4}(\R)$ to obtain the crucial estimates \eqref{notes-2} (ii)
( see remark (b) below). 

 Our first result is concerned with the propagation of regularity in the right hand side of the data for positive times. It affirms that this regularity moves with infinite speed to its left as time evolves. 

\begin{theorem}\label{t1}
If  $u_0\in H^{{3/4}^{+}}(\R)$ and for some $\,l\in \Z^{+},\,\;l\geq 1$ and $x_0\in \R$
\begin{equation}\label{notes-3}
\|\,\partial_x^l u_0\|^2_{L^2((x_0,\infty))}=\int_{x_0}^{\infty}|\partial_x^l u_0(x)|^2dx<\infty,
\end{equation}
then the solution of the IVP \eqref{gkdv} provided by Theorem D satisfies  that for any $v>0$ and $\epsilon>0$
\begin{equation}\label{notes-4}
\underset{0\le t\le T}{\sup}\;\int^{\infty}_{x_0+\epsilon -vt } (\partial_x^j u)^2(x,t)\,dx<c,
\end{equation}
for $j=0,1, \dots, l$ with $c = c(l; \|u_0\|_{{3/4}^{+},2};\|\,\partial_x^l u_0\|_{L^2((x_0,\infty))} ; v; \epsilon; T)$.

In particular, for all $t\in (0,T]$, the restriction of $u(\cdot, t)$ to any interval $(x_0, \infty)$ belongs to $H^l((x_0,\infty))$.

Moreover, for any $v\geq 0$, $\epsilon>0$ and $R>0$ 
\begin{equation}\label{notes-5}
\int_0^T\int_{x_0+\epsilon -vt}^{x_0+R-vt}  (\partial_x^{l+1} u)^2(x,t)\,dx dt< c,
\end{equation}
with  $c = c(l; \|u_0\|_{_{{3/4}^{+},2}};\|\,\partial_x^l u_0\|_{L^2((x_0,\infty))} ; v; \epsilon; R; T)$.
\end{theorem}

 Our second result describes the persistence properties and  regularity effects, for positive times,  in solutions associated with data having polynomial decay in the positive real line.

\begin{theorem}\label{t2}
If $u_0\in H^{{3/4}^{+}}(\R)$ and for some $\,n\in \Z^{+},\;n\geq 1$,
\begin{equation}\label{b1}
\|\,x^{n/2} u_0\|^2_{L^2((0,\infty))}=\int_0^{\infty}| \,x^n|\,|u_0(x)|^2dx<\infty,
\end{equation}
then the solution $u$ of the IVP \eqref{gkdv} provided by Theorem D satisfies
 that
\begin{equation}\label{b2}
\underset{0\le t\le T}{sup}\; \int_0^{\infty} |x^n|\,|u(x,t)|^2\,\,dx\le c
\end{equation}
with $c=c(n;\|u_0\|_{{3/4}^{+},2};\|\,x^{n/2} u_0\|_{L^2((0,\infty))} ; T)$.

Moreover,  for any $\epsilon, \delta, R >0, v\geq 0$, $m,\;j\in\Z^+$, $\;m+j\le n$, $m\ge 1$,
\begin{equation}\label{b3}
\begin{split}
&\underset{\delta \le t\le T}{sup} \;\int_{\epsilon-vt}^{\infty} (\partial_x^mu)^2(x,t)\,x_+^{j}\,dx\\
&\hskip15pt +
\int_{\delta}^T \int_{\epsilon-vt}^{R-vt} (\partial_x^{m+1}u)^2(x,t)\,x_+^{j-1}\,dxdt\le c,
\end{split}
\end{equation}
with $\,c=c(n;\|u_0\|_{{3/4}^{+},2};\|\,x^{n/2} u_0\|_{L^2((0,\infty))} ; T;\delta; \epsilon; R; v)$.
\end{theorem}

It will be clear from our proofs of Theorem \ref{t1} and Theorem \ref{t2} that they still hold for solutions of the \lq\lq defocussing" k-gKdV
\begin{equation*}
\partial_t u+\partial_x^3u-u^k\,\partial_x u=0, \hskip10pt x,t\in\R, \;\;k\in \Z^{+}.
\end{equation*}
Therefore, our results apply to $\,u(-x,-t)\,$ if $\,u(x,t)\,$ is a solution of \eqref{gkdv}. In other words, Theorem \ref{t1} and Theorem \ref{t2} resp. remain valid,
backward in time, for datum satisfying the hypothesis \eqref{notes-3} and \eqref{b1} resp. on the left hand side of the real line.

 As a direct consequence of Theorem \ref{t1} and Theorem \ref{t2}, the above comments   and the time reversible character of the equation in \eqref{gkdv} one has:

\begin{corollary}\label{cor1}
Let $\,u\in C([-T,T] : H^{{3/4}^{+}}(\R))$ be a solution of the equation in \eqref{gkdv} described in Theorem D.
 If there exist $ \,m\in\Z^{+},\;\hat t\in(-T,T),\;a\in\R $ such that
\begin{equation*}
\partial_x^mu(\cdot,\hat t)\notin L^2((a,\infty)),
\end{equation*}
then for any $t\in[-T,\hat t)$ and any $\beta\in\R$ 

\begin{equation*}
%\begin{aligned}
\partial_x^mu(\cdot,t)\notin L^2((\beta,\infty)),
\hskip10pt \text{and} \hskip10pt x^{m/2}\,u(\cdot,t)\notin L^2((0,\infty)).
%\end{aligned}
\end{equation*}
\end{corollary}

 Next, as a consequence of Theorem \ref{t1} and Theorem \ref{t2}  one has  that for an appropriate class of data the singularity of the solution travels with infinite speed to the left as time evolves. In the integrable cases $k=1, 2$ this is expected as part of the so called \it resolution conjecture\rm, (see \cite{EcSc}). Also as a consequence of the time reversibility property one has that the solution cannot have had some regularity in the past.

\begin{corollary}\label{cor2}
Let $\,u\in C([-T,T] : H^{{3/4}^{+}}(\R))$ be a solution of the equation in \eqref{gkdv} described in Theorem D.
 If there exist $ n, m\in\Z^{+}$ with $\,m\leq n\,$ such that for some $a,\,b\in\R$ with $a<b$
\begin{equation}\label{A}
\int_b^{\infty} |\partial_x^nu_0(x)|^2dx <\infty\;\;\;\;\;\;\text{but}\;\;\;\;\;\;\partial_x^mu_0\notin L^2((a,\infty)),
\end{equation}
then for any $t\in(0,T)$ and any $v>0$ and $\epsilon>0$
\begin{equation*}
\int_{b+\epsilon-vt}^{\infty} |\partial_x^nu(x,t)|^2dx <\infty,
\end{equation*}
and for any $t\in(-T,0)$ and any $\alpha\in\R$
\begin{equation*}
\int_{\alpha}^{\infty} |\partial_x^mu(x,t)|^2dx =\infty.
\end{equation*}
\end{corollary}

\underbar{Remark: } If in Corollary \ref{cor2} in addition to  \eqref{A} one assumes that
$$
\int_{-\infty}^a|\partial_x^nu_0(x)|^2dx<\infty,
$$
then by combining the results in this corollary  with the group properties
it follows that
\begin{equation*}
\int_{-\infty}^{\beta} |\partial_x^n u(x,t)|^2\,dx =\infty, \hskip10pt \text{for any}\hskip 5pt \beta\in\R\hskip5pt\text{and}\hskip5pt t>0.
\end{equation*}
This tell us that in general the regularity in the left hand side of the real line does not propagate forward in time.

 The same argument shows  that for initial datum  
 \begin{equation*}
  u_0\in L^2(|x|^ndx)\cap H^{3/4+}(\R)- H^m(\R)
  \end{equation*}
   with $ m, n\in Z^+, \;m<n$, the corresponding solution 
 $\,u(\cdot,t)$ of \eqref{gkdv} satisfies
 \begin{equation*}
x^{n/2}_{+}u(\cdot,t)\in L^2((0,\infty))\;\;\;\;\;\;\text{and}\;\;\;\;\;\;\;x_{-}^{m/2} u(\cdot,t)\notin L^2((-\infty,0)),\;\;\;\;\;\;\;\\t\in(0,T].
 \end{equation*}
 Thus, one has that in general the decay in the left hand side of the real line does not propagate forward in time.

\vskip.1in 
 In a similar manner we also have,

\begin{corollary}\label{cor3}
Let $\,u\in C([-T,T] : H^{{3/4}^{+}}(\R))$ be a solution of the equation in \eqref{gkdv} described in Theorem D.
 If  for $\,j,\,m\in\Z^{+},\;j<m$,
\begin{equation*}
x_{+}^{m/2}\,u_0\in L^2((0,\infty)),\;\;\text{and}\;\; \partial_x^ju_0(x)\notin L^2((\beta,\infty)),\;\;\;\text{for some} \;\; \,\beta\in \R,
\end{equation*}
then for any $t\in(0,T]$
\begin{equation*}
x_{+}^{m/2}\,u(\cdot,t)\in L^2((0,\infty)),\;\;\text{and}\;\; \partial_x^mu(\cdot,t)\in L^2((\alpha,\infty)),\;\;\;\forall \,\alpha\in \R,
\end{equation*}
and for any $\,t\in[-T,0)$
\begin{equation*}
x_{+}^{j/2}\,u(\cdot,t)\notin L^2((0,\infty)),\;\;\text{and}\;\; \partial_x^ju(\cdot,t)\notin L^2((\alpha,\infty)),\;\;\;\forall \,\alpha\in \R.
\end{equation*}
\end{corollary}

 As a consequence of Theorem \ref{t2} we can improve the result in Theorem 1.4 in \cite{ILP1} in the case of a positive integer. More precisely,

\begin{corollary}\label{cor4} Let $u\in C([-T,T]: H^{{3/4}^{+}}(\R))$ be a solution of the equation \eqref{gkdv} described in Theorem D. If there
exist $n_j\in \mathbb{Z^+}\cup\{0\},\;j=1,2,3,4$, $ \,t_0,\,t_1\in[-T,T]$ with $t_0<t_1$ and  $\,a,\,b\in \R$ such that
\begin{equation*}
\int_0^{\infty} \,|x|^{n_1}|u(x,t_0)|^2dx<\infty \hskip10pt and \;\;\;\;\;\;\int_a^{\infty}|\partial_x^{n_2}u(x,t_0)|^2dx<\infty,
\end{equation*}
and
\begin{equation*}
\int_{-\infty}^0 \,|x|^{n_3}|u(x,t_1)|^2dx<\infty \hskip10pt and \;\;\;\;\;\;\;\int_{-\infty}^b|\partial_x^{n_4}u(x,t_1)|^2dx<\infty,
\end{equation*}
then 
$$
u\in C([-T,T]:H^s(\R)\cap L^2(|x|^rdx))
$$
where 
$$
s=\min\{\max\{n_1;n_2\};\max\{n_3;n_4\}\}\;\;\;\;\;\;\text{and}\;\;\;\;\;\;\;\,r=\min\{n_1;n_3\}.
$$
\end{corollary}

\underline {Remarks:} (a)  The improvement of Theorem 1.4 in \cite{ILP1} in the case of a positive follows by taking $\,n_1=n_3\in\Z^+$ and $n_2=n_4=0$.

(b) Although for the sake of the simplicity we shall not pursue this issue here, we remark that the  results in Theorem \ref{t2} and Corollary \ref{cor4} can be extended to non-integer values of the parameter $n$. In this case, one needs to combine the argument given below with those found in \cite{NaPo}.

As it was mentioned above the solution flow of  the IVP \eqref{gkdv} does not preserve the class $Z_{s,r}=H^s(\R)\cap L^2(|x|^rdx)$ when $r>s$. Our next result describes how the decay ($r-s$), not preserved by the flow, is transformed  into extra regularity of the solution. 

\begin{corollary}\label{cor7} If $u_0\in Z_{s,r}$ with $s>3/4^+$, $r\in\Z^+$ and $r>s$. Then  the solution of the IVP  \eqref{gkdv} $u\in C([-T,T]:Z_{s,s})$ also satisfies that for any $b>0$ 
\begin{equation*}
\int_b^{\infty}(\partial_x^ru(x,t))^2dx<\infty,\;\;\;\;\;\;\;\;\text{for any}\;\;t\in(0,T],
\end{equation*}
and
\begin{equation*}
\int_{-\infty}^{-b}(\partial_x^ru(x,t))^2dx<\infty,\;\;\;\;\;\;\;\;\text{for any}\;\;t\in[-T,0).
\end{equation*}

\end{corollary}

\underline {Remarks:} (a) Combining the results in Theorems \ref{t1} and \ref{t2} with those in \cite{ILP1} previously described one can deduce 
further properties of the solution $u(x,t)$ as a consequence of the regularity and decay assumptions on the data $u_0$. For
example, the solution of the IVP \eqref{gkdv} with data $u_0\in H^{{3/4}^{+}}(\R)\backslash H^1(\R)$ for which  there exists
$n\in\Z^{+}$
\begin{equation*}
\int_0^{\infty} |x^n|\,|u_0(x)|^2\,dx<\infty
\end{equation*}
one has, in addition to \eqref{b2}--\eqref{b3}, that for any $M>0$ and any $t>0$
\begin{equation*}
\int_{-\infty}^{-M} |x\,u(x,t)|^2\,dx=\infty \text{\hskip10pt and \hskip5pt} \int_{-\infty}^{-M} |\partial_x\,u(x,t)|^2\,dx=\infty .
\end{equation*}

(b) The proof of Theorems \ref{t1} and \ref{t2} relies on weighted energy estimates and an iterative process for which the estimate \eqref{notes-2} (ii) is essential. This is a consequence of the following version of the Strichartz estimates \cite{Stri} obtained in \cite{KPV89}. The solution 
of the linear IVP 
\begin{equation}
 \begin{aligned}
 \begin{cases}
\label{linearIVP}
 \partial_t v + \partial_x^3 v=0,\\
 v(x,0)=v_0(x),
 \end{cases}
 \end{aligned}
\end{equation}
is given by the group $\{U(t)\,:\,t\in \R\}$
 $$
 U(t)v_0(x)=\frac{1}{\root{3}\of{3t}}\,Ai\left(\frac{\cdot}{\root {3}\of{3t}}\right)\ast v_0(x). 
 $$
where $\,Ai(\cdot)$ denotes the Airy function 
 \begin{equation}
\label{airy}
Ai(x)=c\,\int_{-\infty}^{\infty}\,e^{ ix\xi+i \xi^3/3 }\,d\xi.
 \end{equation}

 The following inequality was established in \cite{KPV89} : for any $(\theta,\alpha)\in [0,1]\times[0,1/2]$ 
\begin{equation}
\label{stric}
\| D^{\theta \alpha/2} U(t)u_0\|_{L^q(\R:L^p(\R))}\leq c\|u_0\|_2,
\end{equation}
with $(q,p)=(6/\theta(1+\alpha),2/(1-\theta))$. In particular, by taking $(\theta,\alpha)=(1,1/2)$ one has
\begin{equation}
\label{stric1}
(\int_{-\infty}^{\infty} \| D^{1/4} U(t)u_0\|_{\infty}^4\,dt)^{1/4}\leq c\|u_0\|_2,
\end{equation}
which explains the hypothesis $s=3/4^+$ in Theorem D and the conclusion \eqref{notes-2} (ii) on it. 

(c) Even for the associated linear IVP \eqref{linearIVP} the results in Theorems \ref{t1} and \ref{t2} do not seem to have been appreciated before, even that in this (linear) case a different proof should follow based on estimates of the Airy function $Ai(\cdot) $ (see \eqref{airy}) and its derivatives.

(d) Although we will not pursue this issue here,  it should be remarked that the results in Theorems \ref{t1} and \ref{t2} can be extended to include some continuity property in time. For example, in Theorem \ref{t1} in addition to \eqref{b2} one can show that for any $t_0\in (0,T)$ and $\epsilon>0$ and any $v>0$ 
$$
\lim_{t\to t_0} \;\int_{x_0+\epsilon-vt}^{\infty}\,(\partial_x^ju)^2(x,t)dx=\int_{x_0+\epsilon-vt_0}^{\infty}\,(\partial_x^ju)^2(x,t_0)dx.
$$
 In this case the proof follows by using an argument similar to that given in \cite{BoSm}.

(e) Without loss of generality from now on we shall assume that in Theorem \ref{t1} $x_0=0$. 

(f) We recall that the above results still hold if one replaces $x, t>0$ by  $x, t<0$. 

(g) In a forthcoming work we shall extend some of the results obtained here to other dispersive models.

The rest of this paper is organized in  the following manner: in section 2 we introduce the family of cut-off functions to be used in the proof of Theorem \ref{t1} and Theorem \ref{t2}. 
This is an important part in our iterative argument used in the proof.
The proof of Theorem \ref{t1} will be given in section 3 and the proof of Theorem \ref{t2} in section 4. 

\section{Preliminaries}

We shall construct  a class of real functions $\car(x)$ for $\epsilon>0$ and $b\ge 5\epsilon$ such that 
\begin{equation*}
\car\in C^{\infty}(\R),\hskip5pt \car'\ge0,
\end{equation*}

\begin{equation*}
\car(x)=\begin{cases}
0,\hskip10pt x\le \epsilon,\\
1, \hskip 10pt x\ge b,
\end{cases}
\end{equation*}
therefore
$$ {\rm supp}\; \car \subseteq [\epsilon,\infty),$$
\begin{equation*}
\car'(x)\geq \frac{1}{b-3\epsilon} \,1_{[3\epsilon, b-2\epsilon]}(x),
\end{equation*}
and
\begin{equation*}
{\rm supp}\; \car'(x)\subseteq [\epsilon, b].
\end{equation*}
Thus 
\begin{equation*}
\chi_{_{0, \epsilon/3, b+\epsilon}}'(x) \ge c_j\, |\car^{(j)}(x)|, \hskip10pt \forall x\in\R,  \;\;\forall j\ge 1.
\end{equation*}
Also, if $\,x\in (3\epsilon, \infty)$, then 
\begin{equation*}
\chi_{_{0, \epsilon, b}}(x)\geq \chi_{_{0, \epsilon, b}}(3\epsilon)\geq \frac{1}{2}\,\frac{\epsilon}{b-3\epsilon},
\end{equation*}
and for any $\,x\in \R$
\begin{equation*}
\chi_{_{0,\epsilon/3,b+\epsilon}}'(x) \leq \frac{1}{b-3\epsilon}.
\end{equation*}

 We also have that given $\,\epsilon>0,\;\;b\geq 5\epsilon$ there exist $c_1,\;c_2>0$ such that
\begin{equation*}
\begin{aligned}
& \chi_{0,\epsilon,b}'(x) \leq c_1\, \chi_{0,\epsilon/3,b+\epsilon}'(x) \chi_{0, \epsilon/3, b+\epsilon}(x),\\
\\
&\chi_{0,\epsilon,b}'(x)   \leq c_2\,\chi_{0, \epsilon/5,\epsilon}(x).
\end{aligned}
\end{equation*}

We shall obtain this family $\{\chi_{_{0, \epsilon, b}}\,:\,\epsilon>0,\;\;b\geq 5\epsilon\}$ by first considering $\rho\in C^{\infty}_0(\R)$, $\rho(x)\geq 0$, even, with $\,\text{supp} \,\rho\subseteq(-1,1)$ and $\,\int\,\rho(x)dx=1$. Then defining
\begin{equation*}
\nu_{_{\epsilon,b}}(x)=\begin{cases}
0,\hskip10pt x\le 2\epsilon,\\
\\
\frac{1}{b-3\epsilon} x-\frac{2\epsilon}{b-3\epsilon},\;\;\;\,x\in[2\epsilon,b-\epsilon],\\
\\
1, \hskip 10pt x\ge b-\epsilon,
\end{cases}
\end{equation*}
and 
\begin{equation*}
\chi_{_{0, \epsilon, b}}(x)=\rho_{\epsilon}\ast \nu_{\epsilon,b}(x),
\end{equation*}
where $\rho_{\epsilon}(x)=\epsilon^{-1}\rho(x/\epsilon)$.

For any $n\in \Z^{+}$ define
\begin{equation*}
\carn(x)= x^n\, \car(x).
\end{equation*}

Thus, for any $n\in \Z^{+}$
\begin{equation*}
\chi_{_{n+1, \epsilon, b}}'(x)\geq (n+1) \chi_{_{n,\epsilon,b}}(x).
\end{equation*}

\section{Proof of Theorem \ref{t1}}
\vspace{3mm}
 We shall use an induction argument. First, we shall prove \eqref{notes-4} for $l=1$ and $l=2$ (to illustrate our method).
 
\vspace{3mm}
\noindent\underline{Case $l=1$}
\vspace{3mm}

Without loss of generality we restrict ourselves to the case $k=1$ in \eqref{gkdv},  (KdV case).

Formally, take partial derivative with respect to $x$ of the equation in \eqref{gkdv} and multiply by $\partial_x u\car(x+vt)$ to obtain 
after integration by parts the identity
\begin{equation}\label{notes-10}
\begin{split}
&\frac12\frac{d}{dt}\int (\partial_x u)^2(x,t)\chi_0(x+vt)\,dx-\underset{A_1}{ \underbrace{v\int  (\partial_x u)^2(x,t)\chi_0'(x+vt)\,dx}}\\
&+\frac32\int  (\partial_x^2 u)^2(x,t)\chi_0'(x+vt)\,dx
-\underset{A_2}{ \underbrace{\frac12\int  (\partial_x u)^2(x,t)\chi_0'''(x+vt)\,dx}}\\
&+ \underset{A_3}{ \underbrace{\int  \partial_x(u\partial_x u)\partial_xu(x,t)\chi_0(x+vt)\,dx}}=0
\end{split}
\end{equation}
where in $\chi_0$ we omit the index $\epsilon, b$ (fixed).

Then after integrating in the time interval $[0,T]$ using the estimate \eqref{notes-2}(iii)  with $r=0$ we have
\begin{equation}\label{notes-11}
\int_0^T |A_1(t)|\,dt \le v\, \int_0^T\int (\partial_x u)^2\chi_0'(x+vt)\,dxdt < c,
\end{equation}
since given $v, \epsilon, b, T$ as above there exist $c_0>0$ and  $R>0$ such that
\begin{equation*}
\chi_0'(x+vt) \le c_0 1_{[-R,R]}(x), \hskip15pt \forall (x,t)\in\R\times [0,T].
\end{equation*}

The same argument shows that
\begin{equation}\label{notes-12}
\int_0^T |A_2(t)|\,dt \le  c_0.
\end{equation}

Finally, since
\begin{equation}\label{notes-13}
\begin{split}
A_3 &=\int \partial_xu \,\partial_xu\,\partial_xu\,\chi_0\,dx+\int u \partial_x^2u\,\partial_xu\,\chi_0\,dx\\
&=\frac12\int  \partial_xu\, \partial_xu\,\partial_xu\,\chi_0(x+vt)\,dx-\frac12 \int u \partial_xu\,\partial_xu\,\chi_0'(x+vt)\,dx\\
&= A_{31}+ A_{32}
\end{split}
\end{equation}
one has 
\begin{equation}\label{notes-14}
|A_{31}|\le \|\partial_x u\|_{\infty}\, \int (\partial_x u)^2 \,\chi_0(x+vt)\,dx
\end{equation}
and
\begin{equation}\label{notes-15}
|A_{32}|\le \|u (t)\|_{\infty}\, \int (\partial_x u)^2 \,\chi_0'(x+vt)\,dx.
\end{equation}

Hence by Sobolev embedding and \eqref{notes-11} after integrating in the time interval $[0,T]$ one gets
\begin{equation}\label{notes-16}
\int_0^T|A_{32}|\,dt \le \sup_{[0,T]} \| u(t)\|_{{3/4}^{+}}\, \int_0^T\int (\partial_x u)^2 \,\chi_0'(x+vt)\,dxdt \le c_0.
\end{equation}

Inserting the above information in \eqref{notes-10}, Gronwall's inequality and \eqref{notes-2} (ii) yield the estimate
\begin{equation}\label{notes-17}
\sup_{[0,T]} \int (\partial_x u)^2\, \car(x+vt)\,dx+\int_0^T\int  (\partial_x^2 u)^2\, \car'(x+vt)\,dxdt\le c_0
\end{equation}
with $c_0= c_0(\epsilon;b;v)>0$  for any $\epsilon>0$, $b\ge 5\epsilon$, $v>0$, which proves the case $l=1$.

\vspace{3mm}
\noindent\underline{Case $l=2$}
\vspace{3mm}

We can assume that the solution $u(\cdot)$ satisfies \eqref{notes-2} and \eqref{notes-17} (case $l=1$). Then formally one
has the following identity
\begin{equation}\label{notes-18}
\begin{split}
&\frac12\frac{d}{dt}\int (\partial_x^2 u)^2(x,t)\chi_0(x+vt)\,dx-\underset{A_1}{ \underbrace{v\int  (\partial_x^2 u)^2(x,t)\chi_0'(x+vt)\,dx}}\\
&+\frac32\int  (\partial_x^3 u)^2(x,t)\chi_0'(x+vt)\,dx
-\underset{A_2}{ \underbrace{\frac12\int  (\partial_x^2 u)^2(x,t)\chi_0'''(x+vt)\,dx}}\\
&+ \underset{A_3}{ \underbrace{\int  \partial_x^2(u\partial_x u)\partial_x^2u(x,t)\,\chi_0(x+vt)\,dx}}=0.
\end{split}
\end{equation}

After integration in time we have from \eqref{notes-17}
\begin{equation}\label{notes-19}
\int_0^T |A_1(t)|\,dt \le |v|\int_0^T\int (\partial_x^2 u)^2\,\chi'_0(x+vt)\,dxdt\le |v|\,c_0.
\end{equation}

Next by the construction of the $\car$'s one has that for $\epsilon>0$, $b\ge 5\epsilon$, there exists $c>0$ such that
\begin{equation}\label{notes-20}
\big| \car'''(x)\big|\le c_{\epsilon, b}\,\chi'_{0,\epsilon/3,b+\epsilon}(x) \hskip15pt \forall x\in \R.
\end{equation}

Therefore, after integration in time using \eqref{notes-17} with $(\epsilon/2, b+\epsilon)$ instead of  $(\epsilon, b)$, it follows that
\begin{equation}\label{notes-21}
\begin{split}
\int_0^T |A_2(t)|\,dt &\le \int_0^T\int (\partial_x^2 u)^2\,|\car'''(x+vt)|\,dxdt\\
&\le c \int_0^T\int (\partial_x^2 u)^2\,|\chi_{0,\epsilon/3, b+\epsilon}'(x+vt)|\,dxdt \le c.
\end{split}
\end{equation}

Finally we have to consider $A_3$ in \eqref{notes-18}. Thus
\begin{equation}\label{notes-22}
\begin{split}
\int &\partial_x^2(u\partial_x u)\,\partial_x^2u\,\chi_0(x+vt)\,dx\\
 &= \int u\,\partial_x^3u\, \partial_x^2 u \chi_0(x+vt)\,dx +3 \int \partial_x u\,\partial_x^2u\, \partial_x^2 u\, \chi_0(x+vt)\,dx\\
&=\frac52  \int \partial_x u\,\partial_x^2u\, \partial_x^2 u \,\chi_0(x+vt)\,dx -\frac12  \int u\,\partial_x^2u\, \partial_x^2 u \,\chi_0'(x+vt)\,dx \\
&= A_{31}+A_{32}.
\end{split}
\end{equation}
Then
\begin{equation}\label{notes-23}
|A_{31}|\le c\|\partial_xu(t)\|_{\infty} \int (\partial_x^2 u)^2\,\chi_0(x+vt)\,dx
\end{equation}
where the last integral is the quantity to be estimated. After integration in time, \eqref{notes-17} and Sobolev embedding we obtain that
\begin{equation}\label{notes-24}
\int_0^T |A_{32}(t)|\,dt \le \sup_{0\le t\le T} \|u(t)\|_{\infty}\int_0^T\int (\partial_x^2 u)^2\,\chi_0'(x+vt)\,dxdt \le c.
\end{equation}

Inserting the above information in \eqref{notes-18} and using Gronwall's inequality it follows that
\begin{equation}\label{notes-25}
\begin{split}
& \sup_{0\le t\le T} \int (\partial_x^2 u)^2\, \car(x+vt)\,dx\\
&+\int_0^T\int  (\partial_x^3 u)^2\, \car'(x+vt)\,dxdt\le c_0
\end{split}
\end{equation}
with $c_0= c_0(\epsilon;b;v)>0$ for any $\epsilon>0$, $b\ge 5\epsilon$, $v>0$.\\

We shall prove the case $l+1$ assuming the case $l\ge 2$. More precisely, we assume:

If $u_0$ satisfies \eqref{notes-3} then \eqref{notes-4} holds, i.e.
\begin{equation}\label{notes-26}
\begin{split}
 &\sup_{0\le t\le T} \int (\partial_x^j u)^2\, \car(x+vt)\,dx\\
&+\int_0^T\int  (\partial_x^{j+1} u)^2\, \car'(x+vt)\,dxdt\le c
\end{split}
\end{equation}
for $j=1,2, \dots, l$, $l\ge2$, for any $\epsilon>0$, $b\ge 5\epsilon$, $v>0$.

Now we have that
\begin{equation}\label{notes-27}
u_0\big|_{(0,\infty)} \in H^{l+1}([0,\infty)).
\end{equation}
Thus from the previous step \eqref{notes-26} holds. Also formally we have for $\epsilon>0$, $b\ge 5\epsilon$, the identity
\begin{equation}\label{notes-28}
\begin{split}
&\frac12\frac{d}{dt}\int (\partial_x^{l+1} u)^2\chi_0(x+vt)dx-\underset{A_1}{ \underbrace{v\int  (\partial_x^{l+1} u)^2\chi_0'(x+vt)dx}}\\
&+\frac32\int  (\partial_x^{l+2} u)^2\chi_0'(x+vt)dx-\underset{A_2}{ \underbrace{\frac12\int  (\partial_x^{l+1}u)^2\chi_0'''(x+vt)dx}}\\
&+ \underset{A_3}{ \underbrace{\int  \partial_x^{l+1}(u\partial_x u)\partial_x^{l+1}u(x,t)\,\chi_0(x+vt)dx}}=0.
\end{split}
\end{equation}

Using \eqref{notes-26} with $j=l$,  after integration in time one has
\begin{equation}\label{notes-29}
\int_0^T |A_1(t)|\,dt\le |v| \int_0^T\int (\partial_x^{l+1} u)^2(x,t)\chi_0'(x+vt)\,dx \le c_0.
\end{equation}

We recall \eqref{notes-20} i.e. given $\epsilon>0$, $b\ge 5\epsilon$, there exists $c>$ such that
\begin{equation}\label{notes-30}
\big|\car'''(x)\big|\le c_{\epsilon, b} \,\chi_{0,\epsilon/3,b+\epsilon}'(x), \hskip15pt \forall x\in\R.
\end{equation}

Hence integrating in the time interval $[0,T]$ and using \eqref{notes-26} with $j=l$ and $(\epsilon, b)=(\epsilon/3,b+\epsilon)$ one
has
\begin{equation}\label{notes-31}
\begin{split}
\int_0^T &|A_2(t)|\,dt \le  \int_0^T\int (\partial_x^{l+1} u)^2(x,t)\,\big|\car'''(x+vt)\big|\,dxdt\\
\le &c\, |v| \int_0^T\int (\partial_x^{l+1} u)^2(x,t)\,\chi_{_{0,\epsilon/3,b+\epsilon}}'(x+vt)\,dxdt \le c_0.
\end{split}
\end{equation}

So it remains to consider $A_3(t)$ in \eqref{notes-28}. Then we have to distinguish  two cases:

\underline{Case $l+1=3$}

Thus 
\begin{equation}\label{notes-32}
\begin{split}
&A_3 = \int \partial_x^3(u\partial_x u)\,\partial_x^3u\, \car(x+vt)\,dx\\
&= 4 \int \partial_x u(\partial_x^3 u)^2\,\car(x+vt)\,dx+ \int u\partial_x^4 u\,\partial_x^3u\, \car(x+vt)\,dx\\
&\hskip10pt +3 \int \partial_x^2u\,\partial_x^2 u\,\partial_x^3u\, \car(x+vt)\,dx\\
&=\frac{7}{2}  \int \partial_x u(\partial_x^3 u)^2\,\car(x+vt)\,dx-\frac12 \int u\,(\partial_x^3u)^2\, \car'(x+vt)\,dx\\
&\hskip10pt + 3\int \partial_x^2u\,\partial_x^2 u\,\partial_x^3u\, \car(x+vt)\,dx\\
&= A_{31}+ A_{32} +A_{33}.
\end{split}
\end{equation}

Then
\begin{equation}\label{notes-33}
|A_{31}(t)|\,dt \le c\|\partial_x u(t)\|_{\infty}\int (\partial_x^3 u)^2\,\car(x+vt)\,dx\\
\end{equation}
where the last integral is the quantity to be estimated. After integration in time and using \eqref{notes-26}
with $j=l=2$ and Sobolev embedding one gets
\begin{equation}\label{notes-34}
\begin{split}
\int_0^T&|A_{32}(t)|\,dt \le c\,\underset{0\le t\le T}{\sup} \|u(t)\|_{\infty}\int_0^T\int (\partial_x^3 u)^2\,\car'(x+vt)\,dxdt\\
&\le \sup_{0\le t\le T}\|u(t)\|_{{3/4}^{+},2}\int_0^T\int (\partial_x^3 u)^2\,\car'(x+vt)\,dx\le c_0.
\end{split}
\end{equation}

So it remains consider $A_{33}(t)$. First we observe that
\begin{equation}\label{notes-35}
\chi_{0,\epsilon/5,\epsilon}(x)=1 \text{\hskip10pt on \hskip5pt} {\rm supp}\; \car\subseteq [\epsilon,\infty)\\
\end{equation}
and since
\begin{equation}\label{notes-36}
\chi_{0,\epsilon/3,b+\epsilon}'(x)\ge c_{\epsilon, b} 1_{[\epsilon, b]}(x) \text{\hskip5pt and \hskip5pt} {\rm supp}\;\car''\subseteq [\epsilon, b]
\end{equation}
one has
\begin{equation}\label{notes-37}
\chi_{0,\epsilon/3,b+\epsilon}'(x)\ge c\big|\car''(x)\big| \hskip15pt \forall x\in\R.
\end{equation}

Thus 
\begin{equation}\label{notes-38}
\begin{split}
A_{33}(t)& =\int \partial_x^2u\,\partial_x^2 u\,\partial_x^3u\, \car(x+vt)\,dx\\
&=-\frac13 \int(\partial_x^2u)^3\, \car'(x+vt)\,dx
\end{split}
\end{equation}
and so by \eqref{notes-35},
\begin{equation}\label{notes-39a}
|A_{33}(t)|\le \|\partial_x^2u\,\car'(\cdot+vt)\|_{\infty} \int (\partial_x^2u)^2\, \chi_{0,\epsilon/5,\epsilon}(x+vt)\,dx
\end{equation}
but the last term is bounded in $ t$ for $t\in [0,T] $ by  a constant $c_{\epsilon, b, v}$ according to \eqref{notes-26} ($j=2$).

Hence
\begin{equation}\label{notes-39b}
\begin{split}
&|A_{33}(t)|\le c\,\|\partial_x^2u\,\car'(\cdot+vt)\|_{\infty} ^2 + c_{\epsilon, b, v}\\
& \le c^{*} \|(\partial_x^2u)^2\,\car'(\cdot+vt)\|_{\infty}  + c_{\epsilon, b, v}, \hskip15pt (c^{*}=c^{*}(\|\car'(\cdot)\|_{L^{\infty}_{x}})\\
&\le c_{\epsilon,b}\int |\partial_x\big((\partial_x^2u)^2\,\car'(x+vt)\big)|\,dx + c_{\epsilon, b, v}\\
&\le c \int |\partial_x^2u\,\partial_x^3u\,\car'(x+vt)|\,dx\\
&\hskip15pt +\int |(\partial_x^2u)^2\,\car''(x+vt)|\,dx + c_{\epsilon, b, v}\\
&\le c \int (\partial_x^2u)^2 \,\car'(x+vt)\,dx+c \int (\partial_x^3u)^2\,\car'(x+vt)\,dx\\
&\hskip15pt + c\int (\partial_x^2u)^2\,\chi_{0,\epsilon/3,b+\epsilon}'(x+vt)\,dx + c_{\epsilon, b, v}.
\end{split}
\end{equation}

By \eqref{notes-26} with $j=1,2$ after integration in time we have that
\begin{equation}\label{notes-40}
\int_0^T |A_{33}(t)|\,dt \le c=c_{\epsilon, b, v,T}.
\end{equation}

\begin{remark}
In the case  $k\ge 2$ in \eqref{gkdv}  when one applies the argument above one basically needs to consider another term of the
form ($k=2$).

\begin{equation}\label{notes-41}
A_4=\int \partial_x u\,\partial_x^2u\,\partial_x^2u\,\partial_x^2u\car(x+vt)\,dx.
\end{equation}

Thus in this case we write (using \eqref{notes-36})
\begin{equation}\label{notes-42}
\begin{split}
&|A_4(t)|\\ 
&\le \int |\partial_xu\,\partial_x^2u|\chi_{0,\epsilon/5,\epsilon}(x+vt)\,dx \;\|(\partial_x^2u)^2\car\|_{\infty}\\
&\le \Big(\!\int \!(\partial_xu)^2\chi_{0,\epsilon/5,\epsilon}(x+vt) dx\!\Big)^{1/2} \Big(\!\int \!(\partial_x^2u)^2\chi_{0,\epsilon/5,\epsilon}(x+vt) dx\!\Big)^{1/2}\\
&\hskip15pt \times \|\partial_x\big((\partial_x^2u)^2\car\big)\|_1\\
&\le  \Big(\int (\partial_xu)^2\chi_{0,\epsilon/5,\epsilon}(x+vt)\,dx + \int (\partial_x^2u)^2\chi_{0,\epsilon/5,\epsilon}(x+vt)\,dx\Big)\\
&\hskip15pt \times\Big( \int \partial_x^2u\,\partial_x^3u\,\car(x+vt)\,dx+ \int (\partial_x^2u)^2\,\car'(x+vt)\,dx \Big).
\end{split}
\end{equation}

Next we observe that
\begin{equation*}
\underset{j=1}{\overset{2}{\sum}} \int (\partial_x^ju)^2\,\chi_{0,\epsilon/5,\epsilon}(x+vt)\,dx
\end{equation*}
is bounded for $ t $ in $[0,T]$ by a constant $c_{\epsilon,b,v, T}=c$.

Coming back to \eqref{notes-42} we have
\begin{equation}\label{notes-42b}
\begin{split}
|A_4(t)| &\le c\big(\int (\partial_x^2u)^2\,\car(x+vt)\,dx+ \int (\partial_x^2u)^2\car'(x+vt)\,dx\\
&\hskip15pt +\int (\partial_x^3u)^2\,\car(x+vt)\,dx\big)\\
&\le c( 1+ \int (\partial_x^2u)^2\car'(x+vt)\,dx\\
&\hskip15pt   +  \int (\partial_x^3u)^2\car(x+vt)\,dx\big)
\end{split}
\end{equation}
where we used that the first integral in the first inequality is  bounded in $[0,T]$. After integration in time the
first integral in the second inequality is bounded by \eqref{notes-26} ($j=1$). Finally, the last integral in the second inequality is the
quantity to be estimated. Gathering the above information yields the desired result.
\end{remark}

\vskip.5mm

\noindent\underline{Case: $l+1\ge 4$ $\;\;(\iff l\ge 3$)}.
\vskip5mm

Returning to \eqref{notes-28} we obtain after integration by parts that
\begin{equation}\label{notes-43}
\begin{split}
A_3&=\int \partial_x^{l+1}(u\partial_xu)\,\partial_x^{l+1}u\,\chi_0(x+vt)\,dx\\
&=c_0 \int u\,\partial_x^{l+1}u\,\partial_x^{l+1}u\,\chi_0'(x+vt)\,dx\\
&\hskip15pt +c_1 \int \partial_x u\,(\partial_x^{l+1}u)^2\,\chi_0(x+vt)\,dx\\
&\hskip15pt +c_2\int \partial_x^2u\, \partial_x^l u\,\partial_x^{l+1}u\,\chi_0(x+vt)\,dx\\
&\hskip15pt +\underset{j=3}{\overset{l-1}{\sum}}\int \partial_x^ju\,\partial_x^{l+2-j}u\,\partial_x^{l+1}u\,\chi_0(x+vt)\,dx\\
&= A_{3,0}+A_{3,1}+A_{3,2}+\underset{j=3}{\overset{l-1}{\sum}} A_{3,j}.
\end{split}
\end{equation}

Next we have
\begin{equation}\label{notes-44}
|A_{3,0}(t)|\le \|u(t)\|_{\infty} \int (\partial_x^{l+1}u)^2\,\chi_0'(x+vt)\,dx
\end{equation}
which after integration in time is bounded by Sobolev embedding and \eqref{notes-26} $j=l$ yield the appropriate bound.

\begin{equation}\label{notes-45}
|A_{3,1}(t)|\le \|\partial_x u(t)\|_{\infty} \int (\partial_x^{l+1}u)^2\,\chi_0(x+vt)\,dx.
\end{equation}
We use the estimate \eqref{notes-2} (ii) to bound the first term on the right hand side of \eqref{notes-45} when we later apply Gronwall's Lemma. The last term is the quantity to be estimated.

To estimate $A_{3,2}$ we follow the argument in the previous case \eqref{notes-38}--\eqref{notes-42}.

So we need to consider  $\underset{j=3}{\overset{l-1}{\sum}} A_{3,j}(t)$ which only appears if $l-1\ge 3$.
Using  an elementary inequality it follows that
\begin{equation}\label{notes-47}
\begin{split}
|A_{3,j}(t)| &= |\int \partial_x^ju\,\partial_x^{l+2-j}u\,\partial_x^{l+1}u \,\car(x+vt)\,dx|\\
&\le \frac12 \int (\partial_x^ju\,\partial_x^{l+2-j}u)^2\,\car(x+vt)\,dx\\
&\hskip15pt +\frac12 \int (\partial_x^{l+1}u)^2 \,\car(x+vt)\,dx\\
&=A_{3,j,1}+  \int (\partial_x^{l+1}u)^2 \,\car(x+vt)\,dx
\end{split}
\end{equation}
where the last integral is the quantity to be estimated. To treat $A_{3,j,1}$ observe that $j$, $l+2-j\leq l-1$ and therefore 
\begin{equation}\label{notes-48}
\begin{split}
|&A_{3,j,1}(t)|\\
&\le \|(\partial_x^{j}u)^2 \chi_{0,\epsilon/5,\epsilon}(\cdot+vt)\|_{\infty} \int (\partial_x^{l+2-j}u)^2\,\car(x+vt)\,dx
\end{split}
\end{equation}
Observe that the last integral is bounded by previous step (induction) \eqref{notes-26}. Also, by Sobolev
\begin{equation}\label{notes-49}
\begin{split}
 &\|(\partial_x^{j}u)^2 \chi_{0,\epsilon/5,\epsilon}(\cdot+vt)\|_{\infty} \le \|\partial_x\Big((\partial_x^{j}u)^2 \chi_{0,\epsilon/5,\epsilon} \Big)\|_1\\
 &\le  \|\partial_x^{j}u\,\partial_x^{j+1}u\,\chi_{0,\epsilon/5,\epsilon}\|_1+  \|(\partial_x^{j}u)^2\, \chi_{0,\epsilon/5,\epsilon}'\|_1\\
 &\le c\Big(\int (\partial_x^{j}u)^2\,\chi_{0,\epsilon/5,\epsilon}(x+vt)\,dx +\int (\partial_x^{j+1}u)^2\,\chi_{0,\epsilon/5,\epsilon}(x+vt)\,dx\\
 &\hskip15pt + \int (\partial_x^{j}u)^2 \chi_{0,\epsilon/5,\epsilon}'(x+vt)\,dx\Big)
 \end{split}
 \end{equation}
 where the first integrals are bounded in the time interval $[0,T]$ (see \eqref{notes-26}) and the last one  is bounded after integration in time  (see \eqref{notes-26}),
 which provides the desired result.

 To justify the previous formal computations we shall follow the following argument. Without loss of generality we assume $l=1$. 
 
 Consider the family of data $u_0^{\mu}=\rho_{\mu} \ast u_0$ with $\rho\in C^{\infty}_0(\R)$, $\text{supp}\,\rho\in (-1,1)$, $\;\rho\ge 0$, $\;\;\displaystyle\int \rho(x)\,dx=1$ and
 \begin{equation*}
  \rho_{\mu}(x)=\frac{1}{\mu} \rho\big(\frac{x}{\mu}\big), \;\;\mu>0.
 \end{equation*} 

For $\mu>0$ consider $u^{\mu}$ the  solutions of the IVP \eqref{gkdv} with data $u_0^{\mu}$  where $(u^{\mu})_{\mu>0}\subseteq C([0,T]: H^{\infty}(\R))$.
 
 Using the continuous dependence of the solution upon the data we have that
\begin{equation}
\label{123}
\sup_{t\in [0,T]}\;\| u^{\mu}(t)-u(t)\|_{3/4^+,2}\,\downarrow 0\;\;\;\;\text{as}\;\;\;\;\mu\,\downarrow 0.
\end{equation}
 Applying the argument in \eqref{notes-10}-\eqref{notes-17} to the smooth solutions $u^{\mu}(\cdot,t)$ one gets that
\begin{equation}\label{124}
\begin{aligned}
\sup_{[0,T]} &\int (\partial_x u^{\mu})^2\, \car(x+vt)\,dx\\
&+\int_0^T\int  (\partial_x^2 u^{\mu})^2\, \car'(x+vt)\,dxdt\le c_0,
\end{aligned}
\end{equation}
for any $\epsilon>0$, $b\ge 5\epsilon$, $v>0$, $c_0= c_0(\epsilon;b;v)>0$  but independent of $\mu>0$
since for $0<\mu<\epsilon$
$$
(\partial_x u_0^{\mu})^2\car(x)=(\partial_x(\rho_{\mu} \ast u_0))^2 \car(x)=(\rho_{\mu}\ast\partial_xu_0\, 1_{[0,\infty)})^2 \car(x).
$$
Combining \eqref{123} and \eqref{124}, and Theorem D (the continuity on the initial data of the class defined in \eqref{notes-2}),  and weak compactness and Fatou's lemma argument one gets that 
\begin{equation}\label{125}
\sup_{[0,T]} \int (\partial_x u)^2\, \car(x+vt)\,dx+\int_0^T\int  (\partial_x^2 u)^2\, \car'(x+vt)\,dxdt\le c_0
\end{equation}
which is the desired result.

\section{Proof of Theorem \ref{t2}}
\vspace{3mm}
\noindent\underline{Proof of \eqref{b2} for any $n\in\Z^{+}$}. First, we observe that 
\vspace{3mm}
\begin{equation}\label{b4}
\begin{split}
x_{+}^n\,u_0\in L^2(\R),&\;\; x_{+}=\max\{x;0\},\;\text{\hskip 5pt implies that for any \hskip5pt}  \epsilon>0  \text{\hskip 5pt and any \hskip5pt}  b\ge 5\,\epsilon, \\
&\hskip20pt \carn(\cdot)\,u_0\in L^2(\R).
\end{split}
\end{equation}

Now using the identity
\begin{equation}\label{b5}
\begin{split}
&\frac12\frac{d}{dt}\int u^2\, \carn(x+vt)\,dx+\frac32\int (\partial_x u)^2 \,\carn'(x+vt)\,dx \\
&+v\, \underset{A_1}{\underbrace{\int u^2\,\carn'(x+vt)\,dx}}\\
&-\frac12  \underset{A_2}{\underbrace{\int u^2 \carn'''(x+vt)\,dx}}+ \underset{A_3}{\underbrace{\int u\partial_x u\, u\,\carn(x+vt)\,dx}}=0.
\end{split}
\end{equation}

We observe that (in general)
\begin{equation}\label{b6}
\begin{split}
\carn'(x)&= (x^n \, \car)'= nx^{n-1}\,\car(\cdot) + x^n\,\car'(\cdot)\\
\carn'''(x)&=(x^n \, \car)'''=n(n-1)(n-2)x^{n-3}\,\car(\cdot) \\
&\hskip10pt +3n(n-1) x^{n-2}\,\car'(\cdot) 
+3n x^{n-1}\,\car''(\cdot)+ x^n\,\car'''(\cdot).
\end{split}
\end{equation}

Thus for any $b>5\epsilon$, there exists $c_{n,b}>0$ such that
\begin{equation}\label{b7}
\big| \carn^{(j)}(x)\big|\le c_{n, b, j} +\carn(x),\hskip10pt  j=1, 2, 3.
\end{equation}

Hence
\begin{equation}\label{b8}
\begin{split}
|A_2(t)| &\le c_{n,b}\int u^2(x,t)\,dx +\int u^2\carn(x+vt)\,dx\\
&\le c_{n,b}\,\|u_0\|_2^2+ \int u^2\carn(x+vt)\,dx.
\end{split}
\end{equation}

Similarly,
\begin{equation}\label{b8b}
\begin{split}
|A_1(t)| &\le |v|\int u^2(x,t)\,\carn(x+vt)\,dx +|v| c_b \int u^2(x,t)\,dx\\
&\le |v|\, \int u^2\carn(x+vt)\,dx+ |v|c_{b}\,\|u_0\|_2^2
\end{split}
\end{equation}
and 
\begin{equation}\label{b9}
|A_3(t)| \le \|\partial_x u(t)\|_{\infty}\, \int u^2(x,t)\,\carn(x+vt)\,dx.
\end{equation}

This together with Gronwall's inequality shows that for any $\epsilon>0$, $\;b\ge 5\epsilon$
\begin{equation}\label{b10}
\begin{split}
&\underset{0\le t\le T}{\sup}\,\int u^2(x,t)\,\carn(x+vt)\,dx\\
&\hskip15pt +\int_0^T\int (\partial_xu)^2(x,t)\,\carn'(x+vt)\,dx <c.
\end{split}
\end{equation}

Next we prove \eqref{b3} in Theorem \ref{t2}.  We divide the proof in several cases. 
\vskip5mm

\noindent\underline{Case $n=1$}

\vskip5mm

From \eqref{b10} and the fact that
\begin{equation}\label{b11}
\chi{_{1,\epsilon,b}}'(x)=(x\,\car)'=\car(x)+ x\,\car'(x) \ge \car(x)
\end{equation}
it follows that for any $\delta>0$ there exists $\hat{t}\in (0,\delta)$ such that
\begin{equation}\label{b12}
\int (\partial_x u)^2(x, \hat{t})\,\car(x)\,dx <\infty.
\end{equation}

Now using the identity
\begin{equation}\label{b13}
\begin{split}
&\frac12\frac{d}{dt}\int (\partial_xu)^2\, \car(x+vt)\,dx+\frac32\int (\partial_x^2 u)^2 \,\car'(x+vt)\,dx \\
&-v\, \underset{A_1}{\underbrace{\int (\partial_xu)^2\,\car'(x+vt)\,dx}}
-\frac12  \underset{A_2}{\underbrace{\int (\partial_x u)^2 \car'''(x+vt)\,dx}}\\
&+ \underset{A_3}{\underbrace{\int\partial_x(u\partial u)\, \partial_xu\,\car(x+vt)\,dx}}=0.
\end{split}
\end{equation}

By \eqref{b10} with $n=1$, after time integration we have that
\begin{equation}\label{b14}
\int_{\hat{t}}^T |A_1(t)|\,dt < c.
\end{equation}

Also using that for any $\epsilon>0$ and $b\ge 5 \epsilon$ there exists $c_{\epsilon,b}$ such that
\begin{equation}\label{b14b}
\big|\car'''(y)\big|\le c_{\epsilon, b}\,\chi_{0,\epsilon/3,b+\epsilon}'(y), \text{\hskip 5pt for all \hskip 5pt} y\in\R,
\end{equation}
from \eqref{b10} with $n=1$ and $(\epsilon/3, b+\epsilon)$ instead of  $(\epsilon, b)$  and after integrating in time one has
\begin{equation}\label{b15}
\int_{\hat{t}}^T |A_2(t)|\,dt < c_{\epsilon, b}.
\end{equation}

Finally, we have (after integration by parts)
\begin{equation}\label{b16}
\begin{split}
\int \partial_x(u\partial_xu)&\,\partial_x u \, \car \,dx \\
&=\frac12 \int \partial_xu\,\partial_xu\,\partial_x u \, \car \,dx -\frac12\int u\partial_x u\,\partial_x u \, \car' \,dx \\
&=A_{3,1}+ A_{3,2},
\end{split}
\end{equation}
where 
\begin{equation}\label{b17}
|A_{3,1}(t)|\le \|\partial_x u(t)\|_{\infty} \int (\partial_x u)^2\,\car(x+vt)\,dx
\end{equation}
and (after integration in time)
\begin{equation}\label{b18}
\begin{split}
\int_{\hat{t}}^T|A_{3,2}&(t)|\, dt \le \int_{\hat{t}}^T \|u(t)\|_{\infty} \int (\partial_x u)^2\,\car'(x+vt)\,dxdt\\
&\le \underset{0\le t \le T}{\sup}\|u(t)\|_{{3/4}^{+}}   \int_{\hat{t}}^T \int (\partial_x u)^2\,\car'(x+vt)\,dxdt < c.
\end{split}
\end{equation}
by \eqref{b10} and Sobolev embedding. Hence from \eqref{b13} and Gronwall's inequality we have that
\begin{equation}\label{b19}
\begin{split}
&\underset{\hat{t}\le t\le T}{\sup}\,\int(\partial_x u)^2(x,t)\,\car(x+vt)\,dx\\
&\hskip15pt +\int_{\hat{t}}^T\int (\partial_x^2u)^2(x,t)\,\car'(x+vt)\,dx <c.
\end{split}
\end{equation}

Before proving the general induction case we shall prove the case $n=2$.

\vskip5mm

\noindent{\underline{Case $n=2$}

\vskip5mm

By hypotheses since $x_{+}\,u_0\in L^2(\R)$ by step 0 we have
\begin{equation}\label{b20}
\begin{split}
&\underset{0\le t\le T}{\sup}\,\int u^2(x,t)\,\chi_{_{2,\epsilon, b}}(x+vt)\,dx\\
&\hskip15pt+\int_{0}^T\int (\partial_x u)^2(x,t)\,\chi_{_{2,\epsilon, b}}'(x+vt)\,dx <c.
\end{split}
\end{equation}

Since
\begin{equation}\label{b21}
\carn'(x)= (x^n\,\car)'= n x^{n-1}\car + x^n \,\car'\ge n \chi_{_{n-1,\epsilon, b}}(x)
\end{equation}
from \eqref{b20} and \eqref{b21} ($n=2$) for any $\delta>0$ there exists $\hat{t}\in(0,\delta)$ such that  $(v\equiv 0$)
\begin{equation}\label{b22}
\int (\partial_xu)^2(x,\hat{t})\,\chi_{_{1,\epsilon, b}}(x)\,dx <\infty.
\end{equation}

Now consider the identity
\begin{equation}\label{b23}
\begin{split}
&\frac12\frac{d}{dt}\int (\partial_xu)^2\, \chi_{_{1}}(x+vt)\,dx+\frac32\int (\partial_x^2 u)^2 \,\chi_{_{1}}'(x+vt)\,dx \\
&-v\, \underset{A_1}{\underbrace{\int (\partial_xu)^2\,\chi_{_{1}}'(x+vt)\,dx}}
-\frac12  \underset{A_2}{\underbrace{\int (\partial_x u)^2 \chi_{_{1}}'''(x+vt)\,dx}}\\
&+ \underset{A_3}{\underbrace{\int\partial_x(u\partial_x u)\, \partial_xu\,\chi_{_{1}}(x+vt)\,dx}}=0.
\end{split}
\end{equation}

By the case $n=1$ using \eqref{b10} after integrating in time we have
\begin{equation}\label{b24}
\int_{\hat{t}}^T |A_1(t)|\,dt \le |v| \int_{\hat{t}}^T\int (\partial_x u)^2(x,t)\,\chi_{_{1,\epsilon,b}}'(x+vt)\,dxdt\le c_v.
\end{equation}

Now since
\begin{equation}\label{b25}
(\chi_1)'''=(x\,\chi_0)'''= 3 \chi_0''+ x\chi_0'''
\end{equation}
we have that given $\epsilon, b, v, T$ there exist $c>0$ and $R>0$ such that
\begin{equation}\label{b26}
\begin{split}
(\chi_1 (x+vt))'''&\le  3 \chi_0''(x+vt)+ x\chi_0'''(x+vt)\\
& \le c\,1_{[-R, R]}(x) \text{\hskip5pt for all \hskip5pt} (x,t)\in \R\times [0,T].
\end{split}
\end{equation}

Thus by \eqref{notes-2} (iii)
\begin{equation}\label{b27}
\int_{\hat{t}}^T |A_2(t)|\,dt \le c \int_{\hat{t}}^T \int_{-R}^R (\partial_x u)^2(x,t)\,dx dt <c.
\end{equation}

%Finally, after integration by parts we obtain
%\begin{equation}\label{b28}
%\begin{split}
%A_3&=\frac12\int \partial_xu\,\partial_x u\, \partial_xu\,\chi_{_{1}}(x+vt)\,dx-\frac12\int u\partial_x u\, \partial_xu\,\chi_{_{1}}'(x+vt)\,dx\\
%&\le \|\partial_x u(t)\|_{\infty} \int (\partial_x u)^2(x,t)\,\chi_{_{1}}(x+vt)\,dx\\
%&\hskip15pt+\|u(t)\|_{\infty}\int (\partial_x u)^2(x,t)\,\chi_{_{1}}'(x+vt)\,dx\\
%&= A_{3,1}+A_{3,2}.
%\end{split}
%\end{equation}
%Then after integration in time it yields
%\begin{equation}\label{b28b}
%\begin{split}
%\int_{\hat{t}}^T& |A_{3,2}(t)|\,dt \\
%&\le \underset{\hat{t}\le t\le T}{\sup}\|u(t)\|_{\infty}\int_{\hat{t}}^T\int (\partial_x u)^2(x,t)\,\chi_{_{1}}'(x+vt)\,dxdt \le c,
%\end{split}
%\end{equation}
%where we used Sobolev embedding and the fact that the last integral is bounded by \eqref{b10} with $n=1$ in the previous step.

To treat $A_3$ we proceed as we did to obtain \eqref{b16} to \eqref{b18}, with $\chi_1$ instead of $\chi_0$. Thus collecting the above 
information we get that
\begin{equation}\label{b29}
\begin{split}
&\underset{\hat{t}\le t\le T}{\sup} \int (\partial_x u)^2(x,t)\,\chi_{_{1,\epsilon,b}}(x+vt)\,dx\\
&\hskip15pt+\int_{\hat{t}}^T\int (\partial_x^2u)^2(x,t)\,\chi_{_{1,\epsilon,b}}'(x+vt)\,dxdt <\infty
\end{split}
\end{equation}
for any $\epsilon>0$ and $b\ge 5\epsilon$.

Now by \eqref{b29} for any $\delta>0$ there exists $\hat{\hat{t}}\in (\hat{t}, \delta)$ such that
\begin{equation}\label{b30}
\int (\partial_x^2 u)^2(x,\hat{\hat{t}}\,)\,\chi_{_{1,\epsilon,b}}'(x)\,dx  <\infty \hskip15pt v\equiv 0,
\end{equation}
which implies 
\begin{equation}\label{b31}
\int (\partial_x^2u)^2(x,\hat{\hat{t}})\,\car(x)\,dx  <\infty .
\end{equation}

Thus
\begin{equation}\label{b32}
\int (\partial^2_x u)^2(x,\hat{\hat{t}}) \car(x)\,dx <\infty, \hskip15pt \forall \epsilon>0, \;\;\forall b\ge 5\epsilon.
\end{equation}

Hence using our result in propagation of regularity (see \eqref{notes-25}) one has: for any $\epsilon>0$, $\,b\ge 5\epsilon$, $v>0$
\begin{equation}\label{b33}
\begin{split}
&\underset{\delta\le t\le T}{\sup} \int (\partial_x^2u)^2(x,t)\,\car(x+vt)\,dx\\
&\hskip15pt+\int_{\delta}^T  \int (\partial_x^3u)^2(x,t)\,\car'(x+vt)\,dx <c.
\end{split}
\end{equation}
This completes the proof of the case $n=2$.

So from now on we shall assume that $n\ge 3$.

To prove Theorem \ref{t2} for the general case $n\ge 3$ we shall use induction.

Given  $(m, l)\in \Z^{+}\times \Z^{+}$ we say that 
\begin{equation}\label{b34}
(m, l) > (\hat{m}, \hat{l}) \iff 
\begin{cases}
{\rm (i)} \hskip10pt m>\hat{m}\\
{\rm or}\\
{\rm(ii)} \hskip10pt m=\hat{m} \text{\hskip5pt and \hskip5pt} l>\hat{l}.
\end{cases}
\end{equation}
Similarly, we say that $ \,(m, l) \geq (\hat{m}, \hat{l})$ if  (ii)  in  the right hand side of \eqref{b34} holds with $\,\geq$ insted of $\,>$.
\vskip5mm
\noindent\underline{Our statement $(m,l)$}
\vskip5mm

For any $\delta>0$, $v>0$, $\epsilon>0$, $b\ge 5\epsilon$
\begin{equation}\label{b35}
\begin{split}
&\underset{\delta\le t\le T}{\sup} \int (\partial_x^l u)^2(x,t)\, \chi_{_{m,\epsilon,b}} (x+vt)\,dx\\
&+ \int_{\delta}^T  \int (\partial_x^{l+1} u)^2(x,t)\, \chi_{_{m,\epsilon,b}}'(x+vt)\,dx \le c_{\delta, v, \epsilon, b}(u_0).
\end{split}
\end{equation}

Notice that we have already  proved the following cases:
\begin{enumerate}
\item $(0,1)$ and $(1,0)$.
\item  $(0,2)$, $(1,1)$ and $(2,0)$.
\item Under the hypothesis $x^{n/2}_{+}\,u_0\in L^2(\R)$ we establish \eqref{b10} i.e. $(n, 0)$ for all $n\in\Z^{+}.$
\item By Theorem A propagation of regularity: If \eqref{b35} holds with $(m,l)=(1,l)$ ($\delta/2$ instead of $\delta$), then there exists $\hat{t}\in (\delta/2,\delta)$ such that
\begin{equation*}
\int (\partial_x^{l+1} u)^2(x, \hat{t})\, \chi_{_{1,\epsilon,b}}'(x+v\hat{t})\,dx <\infty, \hskip10pt v\equiv 0
\end{equation*}
which implies that
\begin{equation*}
\int (\partial_x^{l+1} u)^2(x, t)\, \chi_{_{0,\epsilon,b}}(x+v t)\,dx <\infty \hskip10pt v\equiv 0.
\end{equation*}
By the propagation of regularity (Theorem \ref{t1}) one has the result \eqref{b35} with $(m,l)=(0,l+1)$, i.e.
$(1,l)$ implies $(0,l+1)$ for any $l\in \Z^{+}$.
\end{enumerate}

Hence to complete our induction we assume \eqref{b35} for $(m,k)$ such that 
\begin{equation}\label{b36}
\begin{cases}
{\rm (a)} \hskip 3pt (m,k)\le (n-j, j) \text{\hskip4pt for some \hskip2pt} j=0, 1, \dots, n\\
\text{\hskip20pt and}\\
{\rm (b)} \hskip 3pt (m,k)=(n+1, 0), (n,1), \dots, (n+1-l, l)\hskip2pt\text{for some}\hskip2pt  l\le n.
\end{cases}
\end{equation}

We need to prove the case  $(n+1-(l+1), l+1)=(n-l, l+1)$. Observe that from (4) above, since $(1,l)$ implies $(0,l+1)$, this case is ready for $l=n$. Thus it suffices to
consider $l\le n-1$.

From the previous step ($(n-l+1, l)$) we have that for any $\delta'>0$, $v>0$, $\epsilon$, and $b\ge 5\epsilon$
\begin{equation}\label{b37}
\begin{split}
\underset{\delta'\le t\le T}{\sup}&\, \int (\partial_x^{l} u)^2(x, t)\, \chi_{_{n+1-l,\epsilon,b}}(x+vt)\,dx\\
&+ \int_{\delta'}^T\int (\partial_x^{l+1} u)^2(x, t)\, \chi_{_{n+1-l,\epsilon,b}}'(x+vt)\,dx <\infty.
\end{split}
\end{equation}

Since
\begin{equation}\label{b38}
\chi_{_{n+1-l,\epsilon,b}}'(x)\ge c \, \chi_{_{n-l,\epsilon,b}}(x)
\end{equation}
one has (by \eqref{b37}) that there exists $\hat{t}\in (\delta', 2\delta')$ such that 
\begin{equation}\label{b39}
\int (\partial_x^{l+1} u)^2(x, \hat{t})\, \chi_{_{n-l,\epsilon,b}}(x)\,dx <\infty, \hskip10pt (v\equiv 0).
\end{equation}

Consider the identity  (for $t>\hat{t}$)
\begin{equation}\label{b40}
\begin{split}
\frac12\frac{d}{dt}&\int (\partial_x^{l+1}u)^2(x,t)\, \chi_{_{n-l,\epsilon, b}}(x+vt)\,dx\\
&+\frac32\int (\partial_x^{l+2}u)^2(x,t)\, \chi_{_{n-l,\epsilon, b}}'(x+vt)\,dx\\
&-\underset{A_1}{\underbrace{v\,\int (\partial_x^{l+1}u)^2(x,t)\, \chi_{_{n-l,\epsilon, b}}'(x+vt)\,dx}}\\
&-\underset{A_2}{\underbrace{\frac12\int (\partial_x^{l+1}u)^2(x,t)\, \chi_{_{n-l,\epsilon, b}}'''(x+vt)\,dx}}\\
&+\underset{A_3}{\underbrace{\int \partial_x^{l+1}(u\partial_x u)\, \partial_x^{l+1}u\,\chi_{_{n-l,\epsilon, b}}(x+vt)\,dx}}=0.
\end{split}
\end{equation}

From the previous step $(n-l, l)$ after integration in time $[\hat{t}, T]$ one gets
\begin{equation}\label{b41}
\int_{\hat{t}}^T |A_1(t)|\,dt \le |v|\,\int_{\hat{t}}^T \int (\partial_x^{l+1}u)^2\, \chi_{_{n-l,\epsilon, b}}'(x+vt)\,dxdt \le |v|\,c.
\end{equation}

Next using that
\begin{equation}\label{b42}
\begin{split}
| \chi_{_{n-l,\epsilon, b}}'''(y)|&=|(y^{n-l} \chi_{_{0,\epsilon, b}})'''(y)|\\
&\le c_{nl}  \chi_{_{n-l-3,\epsilon, b}}(y)+\big |c_{nl}^1\, y^{n-l-2}\, \chi_{_{0,\epsilon, b}}'(y)\\
&\;\;\;\;+ c_{nl}^2\, y^{n-l-1} \chi_{_{0,\epsilon, b}}''(y)+ c_{nl}^3 \,y^{n-l} \chi_{_{0,\epsilon, b}}'''(y)\big|\\
&\le c_{nl} \, \chi_{_{n-l-3,\epsilon, b}}(y)+ c_{n,l,b}\,\chi_{_{0,\epsilon/5, \epsilon}}(y).
\end{split}
\end{equation}

Hence
\begin{equation}\label{b43}
\begin{split}
|A_2(t)| &\le c_{ln} \int (\partial_x^{l+1} u)^2\, \chi_{_{n-l-3,\epsilon,b}}(x+vt)\,dx\\
&\;\;\;+
c_{l,n,b} \int (\partial_x^{l+1} u)^2\, \chi_{_{0,\epsilon/5,\epsilon}}(x+vt)\,dx
\end{split}
\end{equation}
which are bounded in $[\hat{t}, T]$ by the step $(n-l-3, l+1)$ and $(0,l+1)$ which is implied by the step $(1,l)= (l+1-l, l)\le (n-l,l)$,
a previous case since $l\le n-1$.

Finally, we need to consider $A_3$ in \eqref{b40}.  Since
\begin{equation}\label{b44}
\begin{split}
\partial_x^{l+1}(u\partial_x u) &= u\,\partial_x^{l+2}u +c_1 \partial_x u \,\partial_x^{l+1}u+ c_2\,\partial_x^2 u\, \partial_x^l u\\
&\hskip15pt +\underset{j=3}{\overset{l-1}{\sum}} c_j \, \partial_x^j u \partial_x^{l+2-j}u.
\end{split}
\end{equation}
we have the following terms in $A_3$
\begin{equation}\label{b45}
\begin{split}
A_3(t)&= \int u\partial_x^{l+2}u \partial_x^{l+1} u\, \chi_{_{n-l,\epsilon, b}}\,dx\\
&\hskip15pt+c_1 \int \partial_x u \,\partial_x^{l+1}u\,\partial_x^{l+1} u\, \chi_{_{n-l,\epsilon, b}}\,dx\\
&\hskip15pt+ c_2\,\int \partial_x^2 u\, \partial_x^l u \partial_x^{l+1} u\,  \chi_{_{n-l,\epsilon, b}}\,dx\\
&\hskip15pt +\underset{j=3}{\overset{l-1}{\sum}} \,c_j \,\int \partial_x^j u\, \partial_x^{l+2-j}u\,\partial_x^{l+1} u\, \chi_{_{n-l,\epsilon, b}}\,dx\\
&= A_{30}+ A_{31} +A_{32} +\underset{j=3}{\overset{l-1}{\sum}}  \; A_{3j} .
\end{split}
\end{equation}

Thus after integration by parts
\begin{equation}\label{b46}
\begin{split}
A_{30}&=-\frac12 \int \partial_x u \,(\partial^{l+1}_x u)^2\, \chi_{_{n-l,\epsilon, b}}\,dx
 -\frac12 \int u \,(\partial^{l+1}_x u)^2\, \chi_{_{n-l,\epsilon, b}}'\,dx\\
& =A_{301}+A_{302} 
\end{split}
\end{equation}
with $A_{301}=-\frac12 A_{31}$ and after integrating in time 
\begin{equation}\label{b47}
\begin{split} 
\int_{\hat{t}}^T |A_{302}|\,dt \le \underset{\hat{t}\le t\le T}{\sup}\, \|u(t)\|_{\infty}\, \int_{\hat{t}}^T \int (\partial_x^{l+1}u)^2 \chi_{_{n-l,\epsilon, b}}'(x+vt)\,dx.
\end{split}
\end{equation}
Notice that the last integral on the right hand side corresponds to the case $(n-l, l)$ (see \eqref{b36}) which is part of our hypothesis of induction and the first term can be bounded as usual using Sobolev embedding.

Next
\begin{equation}\label{b48}
|A_{31}(t)|\le c\,  \|\partial_x u(t)\|_{\infty}\, \int (\partial_x^{l+1}u)^2\, \chi_{_{n-l,\epsilon, b}}(x+vt)\,dx.
\end{equation}
The last term in the inequality above is the quantity to be estimated.

Consider now $A_{32}$ which appears only if $l \ge 2$ (we recall that $n\ge 3$ (to be proved $(n-l, l+1)$))
\begin{equation}\label{49}
A_{32}(t)= c\, \int \partial_x^2u\,\partial_x^lu\,\partial_x^{l+1}u\, \chi_{_{n-l,\epsilon, b}}'(x+vt)\,dx.
\end{equation}
We study two cases:
\vskip5mm
\noindent\underline{ Case $l=2$}
\vskip5mm

By integration by parts (similar to \eqref{notes-39a}--\eqref{notes-40} in the proof of Theorem \ref{t1}) we have
\begin{equation}\label{b50}
\begin{split}
|A_{32}(t)|&=\Big| -\frac{c}3 \int  \int \partial_x^2u\,\partial_x^2u\,\partial_x^2u\, \chi_{_{n-l,\epsilon, b}}'(x+vt)\,dx\Big|\\
&\le c\|(\partial_x^2 u)\chi_{_{n-l,\epsilon, b}}'(\cdot+vt)\|_{\infty} \int (\partial_x^2u)^2\,\chi_{_{0,\epsilon/5, \epsilon}}(x+vt)\,dx\\
&\le c\|(\partial_x^2 u)\chi_{_{n-l,\epsilon, b}}'(\cdot+vt)\|_{\infty}^2+ c\\
&\le c\,c^{*} \,\|(\partial_x^2 u)^2\,\chi_{_{n-l,\epsilon, b}}'(\cdot+vt)\|_{\infty} + c \hskip10pt (c^{*}=\|\chi_{_{n-l,\epsilon, b}}'(\cdot)\|_{L^{\infty}_{x}})\\
&\le c\, \int \big |\partial_x ((\partial_x^2 u)^2\,\chi_{_{n-l,\epsilon, b}}'(x+vt))\big|\,dx +c\\
&\le c\int (\partial_x^2u)^2\, \chi_{_{n-l,\epsilon, b}}'(x+vt)\,dx \\
&\hskip10pt +c\, \int (\partial_x^3u)^2\,\chi_{_{n-l,\epsilon, b}}'(x+vt)\,dx\\
&\hskip10pt+c\,\int (\partial_x^2 u)^2\, \chi_{_{n-l,\epsilon, b}}''(x+vt)\,dx + c\\
&= A_{_{321}}+ A_{_{322}} + A_{_{323}} + c.
\end{split}
\end{equation}
We have used in the second inequality that the second term on the right hand side is bounded in $t\in [\hat{t},T]$ step $(0,l)=(0,2)$.
After integrating in time $A_{321}$ corresponds to the case $(n-l,1)= (n-2,1)$, $A_{322}$ corresponds to the case   $(n-l, 2)= (n-2,2)$ since
$l=2$. So by the hypothesis of induction they are bounded (the previous step to $(n-2, 3)$  which we want to prove).

Finally,
\begin{equation}\label{b51}
\begin{split}
\chi_{_{n-l,\epsilon, b}}''(y)&=(y^{n-l} \,\chi_{_{0,\epsilon, b}})''\\
&=c_{nl} \, y^{n-l-2}\,\chi_{_{0,\epsilon, b}} +c_{nl}^1\, y^{n-l-1}\,\chi_{_{0,\epsilon, b}}' +y^{n-l}\,\chi_{_{0,\epsilon, b}}''.
\end{split}
\end{equation}
Hence there exists $c^{\ast}=c^{\ast}(b, v, T)$ such that
\begin{equation*}
|\chi_{_{n-l,\epsilon, b}}''(x+vt)|\le c\, \chi_{_{n-l-2,\epsilon, b}}(x+vt)+ c_b^{\ast}\, \chi_{_{0,\epsilon/5, \epsilon}}'(x+vt).
\end{equation*}
Thus
\begin{equation}\label{b52}
\begin{split}
|A_{_{323}}(t)| &\le c\int (\partial_x^2u)^2(x,t)\, \chi_{_{n-l-2, \epsilon,  b}}(x+vt)\,dx\\
&\hskip15pt + c_b^{\ast} \int (\partial_x^2u)^2(x,t)\, \chi_{_{0,\epsilon/3, b+\epsilon}},(x+vt)\,dx\\
&= A_{_{3231}}+A_{_{3232}}.
\end{split}
\end{equation}

The term $A_{3231}$ is  bounded in time $t\in[\hat{t},T]$ (case $(n-l-2,2)$) and $A_{3232}$ is also  bounded 
in time (case $(0,2)$).

This handles all the terms in \eqref{b40} for the case $l=2$. Thus after applying Gronwall's inequality one gets
\begin{equation}\label{b53}
\begin{split}
&\underset{\hat{t}\le t\le T}{\sup} \int (\partial_x^{l+1}u)^2(x,t)\, \chi_{_{n-l,\epsilon, b}}(x+vt)\,dx\\
&\hskip5pt +\int_{\hat{t}}^T \int  (\partial_x^{l+2}u)^2(x,t)\, \chi_{_{n-l,\epsilon, b}}'(x+vt)\,dxdt < c.
\end{split}
\end{equation}

It remains to consider the case $l\ge 3$ in \eqref{b45}.
\vskip5mm
\noindent\underline{Case $l\ge 3$}
\vskip5mm

Observe that in this case after integration by parts
\begin{equation}\label{b54}
\begin{split}
A_{32}&= c_2 \int \partial_x^2 u\,\partial_x^{l}u\,\partial_x^{l+1} u\, \chi_{_{n-l,\epsilon, b}}(x+vt)\,dx\\
&=-\frac{c_2}{2}  \int \partial_x^3 u\,(\partial_x^{l}u)^2\, \chi_{_{n-l,\epsilon, b}}(x+vt)\,dx\\\
&\hskip15pt + \frac{c_2}{2} \int \partial_x^2 u\,(\partial_x^{l}u)^2\,\, \chi_{_{n-l,\epsilon, b}}'(x+vt)\,dx
\end{split}
\end{equation}
whose bound is similar to that given above for the case $l=2$. Thus it suffices to consider the reminder terms in \eqref{b45}, i.e.
\begin{equation}\label{b55}
A_{3j}(t)= c_j\,\int \partial_x^ju\,\partial_x^{l+2-j}u\,\partial_x^{l+1}u\,\chi_{_{n-l,\epsilon, b}}(x+vt)\,dx
\end{equation}
without loss of generality (since $l\ge 3$ and $l+2=j+(l+2-j)$ we can assume $3\le j\le l/2+1$.

Thus 
\begin{equation}\label{b56}
\begin{split}
|A_{3j}(t)|&\le c_j \int (\partial_x^ju\,\partial_x^{l+2-j}u)^2\, \chi_{_{n-l,\epsilon, b}}(x+vt)\,dx\\
&\hskip15pt + c_j \int (\partial_x^{l+1}u)^2\, \chi_{_{n-l,\epsilon, b}}(x+vt)\,dx
\end{split}
\end{equation}
where the second term is the quantity to be estimated.

So we need to bound
\begin{equation}\label{b57}
\begin{split}
&\int (\partial_x^ju)^2\,(\partial_x^{l+2-j}u)^2\, \chi_{_{n-l,\epsilon, b}}(x+vt)\,\chi_{_{0,\epsilon/5, \epsilon}}(x+vt)\,dx\\
&\le \| (\partial_x^ju)^2\,\chi_{_{0,\epsilon/5, \epsilon}}(\cdot+vt)\|_{\infty}  \int (\partial_x^{l+2-j}u)^2\, \chi_{_{n-l,\epsilon, b}}(x+vt)\,dx.
\end{split}
\end{equation}
Noticing that the second term is  bounded in time  (case $(n-l,l+2-j)$, $j\ge 3$) by the induction hypothesis.

Thus it remains to bound
\begin{equation}\label{b58}
\begin{split}
& \| (\partial_x^ju)^2\,\chi_{_{0,\epsilon/5, \epsilon}}(\cdot+v t)\|_{\infty}\le \int |\partial_x ((\partial_x^ju)^2\,\chi_{_{0,\epsilon/5, \epsilon}}(\cdot+vt))|\,dx\\
 &\le c\, \int |\partial_x^ju \,\partial_x^{j+1}u\,\chi_{_{0,\epsilon/5, \epsilon}}(\cdot+vt)|\,dx+ \int |(\partial_x^ju)^2\,\chi_{_{0,\epsilon/5, \epsilon}}'(\cdot+vt)|\,dx\\
 &\le \int (\partial_x^ju)^2\,\chi_{_{0,\epsilon/5, \epsilon}}(x+vt)\,dx+\int (\partial_x^{j+1}u)^2\,\chi_{_{0,\epsilon/5, \epsilon}}(x+vt)\,dx\\
 &\hskip15pt +\int (\partial_x^ju)^2\,\chi_{_{0,\epsilon/5, \epsilon}}'(x+vt)\,dx.
 \end{split}
 \end{equation}
 
 The first two terms are bounded in time, cases $(0,j)$, $(0,j+1)$ since $j\le l-1$ implies $j+1\le l$ and $l\le n$.
 
 For the third one after integration in time one gets
 \begin{equation}\label{b59}
 \int_{\hat{t}}^T\int (\partial_x^ju)^2\,\chi_{_{0,\epsilon/5, \epsilon}}'(x+vt)\,dxdt
 \end{equation}
 corresponding to the case $(0,j-1)$ previous case.
 
 This basically completes the proof of Theorem \ref{t2}.

The previous formal computation can be justified by approximating the initial data $u_0$ by elements in the Schwartz space $\mathcal S(\R)$
say $u_0^{\mu}$, $\mu>0$. This is done by an appropriate truncation and smoothing proccess.
Using the well-posedness in this class one gets a family of solutions $u^{\mu}(\cdot,t)$ for which each step of the  above argument can be justified. Due to our construction the 
estimates are uniform in the parameter $\mu>0$ which yields the desired estimate by passing to the limit.

\section*{Acknowledgments}
P. I. was supported by Universidad Na\-cio\-nal de Co\-lom\-bia-Me\-de\-ll\'in. F. L. was partially supported
by CNPq and FAPERJ/Brazil. G. P. was  supported by a NSF grant  DMS-1101499.  Part of this work was complete
while the third author was visiting IMPA under the sponsorship  of the program Ci\^encia sem Fronteiras.

\end{document}